\documentclass[11pt]{article}

\usepackage{amsmath,epsfig,epsf,psfrag, url}
\usepackage{amsmath, amsthm, amssymb}
\usepackage{multirow,color}

\evensidemargin \oddsidemargin

\newtheorem{theorem}{Theorem}[section]

\begin{document}
\setlength{\textheight}{575pt} 
\title{Optimal Prediction in an Additive Functional Model}
\author{Xiao Wang\thanks{Xiao Wang is Associate Professor, Department of Statistics, Purdue University, West Lafayette, IN 47907-2066, USA. (E-mail: wangxiao@purdue.edu)} \and David Ruppert\thanks{David Ruppert is Andrew Schultz Jr.\ Professor of Engineering and
Professor of Statistical Science, Department of Statistical Science and School of Operations Research
and Information Engineering, Cornell University, Comstock Hall, Ithaca, NY
14853, USA. (E-mail: dr24@cornell.edu) }}
\maketitle

\begin{abstract}
The functional generalized additive model (FGAM) provides a more flexible nonlinear functional regression model than the well-studied functional linear regression model.
This paper restricts attention to the FGAM with identity link and additive errors, which we will call the additive functional model, a generalization of the functional linear model.
This paper studies the minimax rate of convergence of predictions from the additive functional model  in the framework of reproducing kernel Hilbert space. It is shown that the optimal rate is determined by the decay rate of the eigenvalues of a specific kernel function, which in turn is determined by the reproducing kernel and the joint distribution of any two points in the random predictor function. For the special case of the functional linear model, this kernel function is jointly determined by the covariance function of the predictor function and the reproducing kernel.
The easily implementable roughness-regularized predictor is shown to achieve the optimal rate of convergence. Numerical studies are carried out to illustrate the merits of the predictor. Our simulations and real data examples demonstrate a competitive performance against the existing approach.\\

\noindent{\it Keywords}: Functional regression, minimax rate of convergence, principal component analysis, reproducing kernel Hilbert space.
\end{abstract}

\section{Introduction}

Functional regression, in particular functional linear regression, has been studied extensively. Recent synopses include \cite{ramsay_02, ramsay_05}, \cite{ferraty_06}, and \cite{ramsay_09}.
Let $X(\cdot)$ be a random process defined on $[0, 1]$ and $Y$ be the univariate response variable.  Typically, $t$ is restricted to a compact interval, so the assumption that $t \in [0,1]$ causes no loss of generality.  
Suppose we observe $n$ i.i.d.\ copies of $\big(Y, X\big)$, $\big(Y_i, X_i\big)$, $i=1, \ldots, n$. The functional linear regression model assumes that
\begin{equation}\label{equ:flr}
Y_i = \alpha_0 + \int_0^1 \beta_0(t) X_i(t) dt + \epsilon_i,
\end{equation}
where $\alpha_0\in \mathbb R$ is the coefficient constant, $\beta_0: [0, 1]\rightarrow \mathbb R$ is the slope function, and the $\epsilon_i$ are i.i.d.\ random errors with $\mathbb E\epsilon_i =0$ and $\mathbb E \epsilon_i^2=\sigma^2, \ 0 < \sigma^2 < \infty$.
One of the popular methods for estimating functional linear models is based on functional principal component analysis (see, e.g., \cite{james_02}, \cite{ramsay_05}, \cite{yao_05}, \cite{cai_06}, \cite{li_07}, \cite{hall_07}). In addition,  methods of regularization have also been applied to the functional linear model (see, e.g., \cite{crambes_09}, \cite{yuan_10}, \cite{cai_12}).

Due to the limitation of the inherent linearity of (\ref{equ:flr}), \cite{ferraty_06} extended this model to nonparametric functional models and \cite{muller_08} discussed functional models that are additive in the functional principal component scores of the predictor functions. Recently, \cite{mclean_12}   proposed a new model called a functional generalized additive model (FGAM).  The same model was studied by \cite{muller_12} who called it the continuously additive model.
We will study the special case of the FGAM with the identify link and continuous errors so that
\begin{equation}\label{equ:funadd}
Y_i = \int_0^1 F_0\Big(t, X_i(t)\Big)dt + \epsilon_i,
\end{equation}
where $F_0(\cdot, ~\cdot): [0, 1]^2\rightarrow \mathbb R$ is a bivariate function.  Because $F_0$ is nonlinear, $X(t)$ can be replaced by $G\{X(t)\}$ for a transformation $G$.  Since $G$ can be strictly increasing function from the entire real line to $[0,1]$, assuming that $X(t) \in [0,1]$ also causes no loss of generality.  (In \cite{mclean_12}, $G_t$ is allowed to depend on $t$ and is an estimate of the CDF of $X(t)$, but we will not pursue this refinement here.)  Model \eqref{equ:funadd} will be called the additive functional model and contains (\ref{equ:flr}) as a special case with $F_0(t, x) =\alpha_0 + x \beta_0(t)$. The additive functional model  offers increased flexibility compared to (\ref{equ:flr}), while still facilitating interpretation and estimation.  In \cite{mclean_12},  computational issues of this model were studied and  $F_0$ was estimated using tensor-product B-splines with roughness penalties. In
\cite{muller_12}, a piecewise constant function was fit to $F_0$ and the asymptotic properties, e.g., consistency and asymptotic normality, of  predictions based on $\widehat F_0$ were studied.

In this paper, we study the minimax prediction. The unknown bivariate function $F_0$ is assumed to reside in a RKHS ${\cal H}(K)$ with a reproducing kernel $K: [0,1]^2\times [0,1]^2\rightarrow \mathbb R$.
The goal of prediction is to recover the functional $\eta_0$:
$$\eta_0(X) = \int_0^1 F_0\Big(t, X(t)\Big)dt,$$
based on the training sample $(Y_i, X_i)$, $i=1,\ldots, n$. Let $\widehat F_n$ be an estimate of $F_0$ from the training data. Then its accuracy can be naturally measured by the excess risk:
\begin{align*}
{\mathfrak R}_n :=& \mathbb E^*\Big[Y_{n+1} - \int_0^1 \widehat F_n\Big(t, X_{n+1}(t)\Big)dt\Big]^2 - \mathbb E^*\Big[Y_{n+1} - \int_0^1 F_0\Big(t, X_{n+1}(t)\Big)dt\Big]^2 \\
=&   \mathbb E^*\Big\{ \int_0^1 \Big[ \widehat F_n(t, X_{n+1}(t)) - F_0(t, X_{n+1}(t))   \Big]dt \Big\}^2,
\end{align*}
where $(Y_{n+1}, X_{n+1})$ possesses the same distribution with $(Y_i, X_i)$ and is independent with $(Y_i, X_i)$, $i=1,\ldots, n$, and $\mathbb E^*$ represents taking expectation over $(Y_{n+1}, X_{n+1})$ only. It is interesting to study the rate of convergence of ${\mathfrak R}_n$ as the sample size $n$ increases, which reflects the difficulty of the prediction problem. A closed related but different problem is estimation the bivariate function $F_0$.

The optimal rate of convergence for the prediction problem is established in this paper. The spectral theorem admits that there exist a set of orthonormalized eigenfunctions $\{\psi_k: k\ge 1\}$ and  a sequence of eigenvalues $\kappa_1\ge \kappa_2\ge \cdots>0$ such that
$$
K\Big((t,x); (s,y)\Big) = \sum_{k=1}^\infty \kappa_k \psi_k(t,x)\psi_k(s,y), ~~K(\psi_k) :=\int\int K\Big(\cdot; (s,y)\Big)\psi_k(s,y)dsdy = \kappa_k \psi_k.
$$
It is shown that under  model (\ref{equ:funadd}), the difficulty of the prediction problem as measured by the minimax rate of convergence depends on the decay rate of the eigenvalues of the kernel $C: [0,1]^2 \times [0, 1]^2 \rightarrow \mathbb R$, and
\begin{align}
 C\Big((t,x); (s,y)\Big)   := \int\int \mathbb E \Big\{K^{1/2}\Big((t,x); (u,X(u))\Big)~~K^{1/2}\Big((s,y); (v,X(v))\Big)\Big\}dudv \label{equ:C}
\end{align}
where  $K^{1/2}\Big((t,x); (s,y)\Big) = \sum_{k=1}^\infty \kappa_k^{1/2} \psi_k(t,x)\psi_k(s,y)$. A minimax lower bound is first derived for the prediction problem. Then a roughness-regularized predictor is introduced and is shown to attain the rate of convergence given in the lower bound. Therefore,  this estimator is rate-optimal.

The paper is organized as follows. Section 2 establishes the minimax lower bound for the rate of convergence of the excess risk. Section 3 develops a predictor using a roughness regularization method and shows this predictor is rate-optimal. Section 4 conducts a Monte Carlo study to validate the method and we also illustrate the merit of the method by using two real data examples. Some discussions are provided in Section 5. The paper ends with proofs in Section 6.

\section{Minimax Lower Bound}

In this section, we establish the minimax lower bound for the rate of convergence of the excess risk.

Assume that the unknown $F_0$ resides in a reproducing kernel Hilbert space ${\cal H}(K)$ with a reproducing kernel $K$. It is well-known that ${\cal H}(K)$ is a linear functional space endowed with an inner product $\langle\cdot, \cdot\rangle_{{\cal H}(K)}$ such that
$$F(t,x) = \Big\langle K\Big((t, x); \cdot\Big), F\Big\rangle_{{\cal H}(K)}, ~~~ \mbox{for any } F\in {\cal H}(K).$$
There is a one-to-one relationship between $K$ and ${\cal H}(K)$.
It follows from (\ref{equ:C}) that
\begin{align*}
&C\Big((t,x); (s,y)\Big)\\
 &= \int\int\int\int\ \bigg\{ K^{1/2}\big((t,x); (u,z_1)\big)~K^{1/2}\big((s,y); (v,z_2)\big)  g\big((u, z_1);(v, z_2)\big)\bigg\}dudvdz_1dz_2,
\end{align*}
where $g\big((u, z_1);(v, z_2)\big)$ is the joint density function of $(X(u), X(v))$ evaluated at $(z_1, z_2)$. Similarly, $C$ admits the spectral decomposition,
$$C\big((t,x); (s,y)\big) = \sum_{j=1}^\infty \rho_j \phi_j(t, x) \phi_j(s, y),$$
where the $\rho_j$ are the positive eigenvalues with a decreasing order and the $\phi_j$ are the corresponding orthonormal eigenfunctions.
We assume $\rho_k \asymp k^{-2r}$ for some constant $0<r<\infty$, where for two sequences $a_k, b_k>0$, $a_k\asymp b_k$ means that $a_k/b_k$ is bounded away from zero and infinity as $k\rightarrow\infty$.


\begin{theorem}\label{thm:low}
Suppose that the eigenvalues $\{\rho_k: k\ge 1\}$ of the kernel $C$ in (\ref{equ:C}) satisfy $\rho_k \asymp k^{-2r}$ for some constant $0<r<\infty$, then the excess prediction risk satisfies
\begin{equation}
\lim_{c\rightarrow 0}\lim_{n\rightarrow\infty} \inf_{\tilde \eta} \sup_{F_0\in {\cal H}(K)} \mathbb P\Big( {\mathfrak R}_n \ge c n^{-{2r\over 2r+1}}  \Big) =1,
\end{equation}
where the infimum is taken over all possible predictors $\tilde \eta$ based on $\{(Y_i, X_i): i=1, \ldots, n\}$.
\end{theorem}

It is interesting to compare Theorem \ref{thm:low} with some of the known results when  functional linear regression  is the true model. If the bivariate function $F$ is restricted to the specific form $F(t, x) = \beta(t) x$, where $\beta$ belongs to a reproducing kernel Hilbert space ${\cal H}(\tilde K)$ with the reproducing kernel $\widetilde K: [0,1]\times [0, 1]\rightarrow \mathbb R$, then we have a functional linear regression model.
Assume
$\widetilde K(t, s) = \sum_{k=1}^\infty \varsigma_k \varphi_k(t)\varphi_k(s),$
where the $(\varsigma_k, \varphi_k)$ are the eigenvalue and eigenfunction pairs for $\widetilde K$. It is not hard to see that
$K\Big((t, x); (s, y)\Big) = 3 \widetilde K(t, s) x y = \sum_{k=1}^\infty \kappa_k\psi_k(t, x)\psi_k(s, y),$
where
$\kappa_k  = \varsigma_k, ~~~ \psi_k(t, x) = \sqrt{3} x \varphi_k(t).$
Therefore,
$$C\Big((t,x); (s,y)\Big) = 3 x y \int\int \widetilde K^{1/2}(t, u) G(u, v) \widetilde K^{1/2}(v, s)dudv,$$
where $G(u, v) = \mathrm{cov}(X(u), X(v))$ is the covariance function of $X$, so the eigenvalues of $C$ have the same decay rate as the eigenvalues of $\widetilde K^{1/2} G\widetilde K^{1/2}$. This special setting coincides with those considered in \cite{yuan_10} and \cite{cai_12}. Results similar to ours have been established  in these papers for this special setting.

\section{A Roughness Regularized Estimate}

In this section, we will develop a predictor using a roughness regularization method and establish that this predictor achieves the optimal rate established in Theorem \ref{thm:low}.

\subsection{Computation}

We define the estimate $\widehat F_{n\lambda}$ of $F_0$ as the minimizer of the functional
\begin{equation}\label{equ:obj}
 {1\over n}\sum_{i=1}^n\Big( Y_i - \int_0^1 F(t, X_i(t))dt \Big)^2 + \lambda J(F),
\end{equation}
where $\lambda$ is the tuning parameter and  $J(\cdot)$ is a squared semi-norm on ${\cal H}(K)$. The first term measures the closeness of the fit  to the data, the second term controls the smoothness of the estimate, and the tuning parameter $\lambda$ adjusts the trade-off between these two.
The estimate $\widehat F_{n\lambda}$ can be computed explicitly over the infinitely dimensional function space ${\cal H}(K)$. This observation is important to both numerical implementation of the procedure and our asymptotic analysis.

Let ${\cal H}_0$ be the null space of $J$, i.e., ${\cal H}_0 = \{F\in {\cal H}: J(F)=0\}$. Assume that $\{\xi_1, \ldots, \xi_N\}$ be the orthonormal basis of ${\cal H}_0$ with $N = \mathrm{dim}({\cal H}_0)<\infty$. Let ${\cal H}_1$ be its orthogonal complement in ${\cal H}$ such that ${\cal H} = {\cal H}_0 \oplus {\cal H}_1$.
\begin{theorem}\label{thm:sol}
The minimizer of (\ref{equ:obj}) over ${\cal H}(K)$ can be represented by
\begin{equation}
\widehat F_{n\lambda}(t, x) = \sum_{j=1}^N d_j \xi_j(t, x) + \sum_{i=1}^n c_i \int_0^1 K\Big( (t, x);(s, X_i(s)) \Big)ds,
\end{equation}
for some  $c = (c_1, \ldots, c_n)^T\in \mathbb R^n$ and $d = (d_1, \ldots, d_N)^T\in \mathbb R^N$.
\end{theorem}

Denote by $\Sigma$ the $n\times n$ matrix with
$(\Sigma)_{ij} = \int\int K\Big((t, X_j(t)); (s, X_i(s))\Big)dtds,$
and
by $\Xi$ the $n\times N$ matrix with
$(\Xi)_{ij} = \int \xi_j(t, X_i(t))dt.$
Then, (\ref{equ:obj}) may be written as the matrix form
\begin{equation}\label{equ:obj2}
 {1\over n} \|Y - \Xi d - \Sigma c\|_2^2 + \lambda c^T \Sigma c,
 \end{equation}
where $J(F) =c^T \Sigma c$. 
 It is easy to see that the solution of the linear system
\begin{align}
( \Sigma + n\lambda I ) c + \Xi d =& Y, \label{equ:1}\\
\Xi^T\Sigma c + \Xi^T\Xi d =& \Xi^T Y,\label{equ:2}
\end{align}
is a solution of (\ref{equ:obj2}). It follows from (\ref{equ:1}) and (\ref{equ:2}) that $\Xi^T c=0$.
Suppose $\Xi$ is of full column rank. Let
$$\Xi = Q R^* = (Q_1, Q_2) \left(\begin{array}{cc}R\\ 0\end{array}\right) = Q_1 R$$
be the QR-decomposition of $\Xi$ with $Q$ orthogonal and $R$ upper-triangular. From $\Xi^T c=0$,  $Q_1^T c=0$, so $c \perp \text{row}(Q_1)$, the row space of $Q_1$. Since $Q$ is orthogonal, $c \in \text{row}(Q_2)$, and $c = Q_2Q_2^T c$ because $Q_2 Q_2^T$ projects onto  $\text{row}(Q_2)$. Simple algebra gives
\begin{align*}
\widehat c & = Q_2(Q_2^T \Sigma Q_2 + n\lambda I)^{-1} Q_2^T Y,\\
\widehat d &= R^{-1}(Q_1^T Y - Q_1^T \Sigma c).
\end{align*}

\subsection{Rate of convergence}

In this section, we turn to the asymptotic properties of the estimate $\widehat F_{n\lambda}$.

\begin{theorem}\label{thm:upper} Assume that for any $F\in L_2([0,1]^2)$
\begin{equation}\label{equ:condition}
\mathbb E\Big( \int F(t, X(t))dt \Big)^4 \le c\Big( \mathbb E\Big( \int F(t, X(t))dt \Big)^2  \Big)^2
\end{equation}
for a positive constant $c$. Then,
\begin{equation}
\lim_{A\rightarrow \infty}\lim_{n\rightarrow\infty}\sup_{F_0\in {\cal H}(K)}\mathbb P\Big\{ \mathfrak{R}_n \ge A n^{-{2r\over 2r+1}}  \Big\} = 0,
\end{equation}
when $\lambda$ is of order $n^{-2r/(2r+1)}$.
\end{theorem}

We have made an additional assumption (\ref{equ:condition}) on $X$. For the functional linear regression model when $F(t, x) = \beta(t) x$, condition (\ref{equ:condition}) shows that, for any $\beta\in L_2([0, 1])$,
$\mathbb E\big( \int \beta(t)X(t)dt \big)^4 \le c\Big( \mathbb E\big( \int \beta(t) X(t)dt \big)^2  \big)^2,$
which states that linear functionals of $X$ have bounded kurtosis. In general, (\ref{equ:condition}) states that such special nonlinear functional $F\big(\cdot, X(t)\big)$ of $X$ have bounded kurtosis.

It follows from both Theorem \ref{thm:low} and Theorem \ref{thm:upper} that the minimax rate of convergence for the excess prediction ${\mathfrak R}_n$ is of order $n^{-2r/(2r+1)}$, which is determined by the decay rate of the eigenvalues of the kernel $C$.

\subsection{Optimal choice of $\lambda$}

Let $\widehat Y = \Big( \eta_{\widehat F_\lambda}(X_1), \ldots,   \eta_{\widehat F_\lambda}(X_n)\Big)^T$. Since the regularized estimator is a linear estimator in $Y$, $\widehat Y = H(\lambda) Y$, where $H(\lambda)$ is called the hat matrix depending on $\lambda$. Some algebra yields
$$H(\lambda) = I - n\lambda F_2 (F_2^T \Sigma F_2+n \lambda I)^{-1} F_2^T.$$
We may select the tuning parameter $\lambda$ that minimizes the generalized cross-validation score \cite{wahba_90},
\begin{equation}\label{equ:gcv}
\mbox{GCV}(\lambda) = {\|\widehat Y - Y\|_2^2/n \over \Big\{1-\mbox{tr}(H(\lambda))/n\Big\}^2}.
\end{equation}
Choosing  $\lambda$ by minimizing GCV worked very well in our numerical studies.

\section{Numerical Results}

In our numerical studies, we compare the numerical performance of the proposed predictor with some well-known existing predictors.

We will focus on a RKHS ${\cal H}(K)$ with a squared seminorm
$$
J(F) = \sum_{\alpha_1+\alpha_2=m} {m!\over \alpha_1!\alpha_2!}\int\int \Big({\partial^m F\over \partial t^{\alpha_1}\partial x^{\alpha_2}}\Big)^2dtdx.
$$
The function
$J_m\Big((t-x)^2 + (x-y)^2\Big),$
where $J_m(x) = x^{2m-2}\log x$ acts like a reproducing kernel in this approach to the computation of thin-plate splines, and hence is called a semi-kernel (\cite{duchon_77}, \cite{meinguet_79}).
In this setting,
the optimal solution of the roughness-regularized estimate can be written as
\begin{align}
F(t, x) = \sum_{j=1}^N d_j \xi_j(t, x) + \sum_{i=1}^n c_i \int J_m\Big( \sqrt{ (t-s)^2 +(x-X_i(s))^2}~  \Big)ds,
\end{align}
where $\xi_j(t, x) = t^{\gamma_1} x^{\gamma_2}$ for some pair of integers $\gamma_1,\gamma_2$ with $0\le \gamma_1+\gamma_2<m$ and $N$ is the number of such pairs. Let $\hat c$ and $\hat d$ be the estimates from the training data. Then, for any random function $X$, the predicted response is
$$\eta_{\widehat F}(X) = \sum_{j=1}^N \hat d_j \int \xi_j(t, X(t))dt + \sum_{i=1}^n \hat c_i \int\int J_m\Big( \sqrt{ (t-s)^2 +(X(t)-X_i(s))^2}~  \Big)dtds.$$
In particular, when $m=2$, we have $N=3$, and
$$\xi_1(t, x) =1, ~~~\xi_2(t, x) =t,~~~ \xi_3(t, x) =x,~~~J_m(x) = x^{2}\log x.$$
Note that $\int \xi_1(t, X(t))dt =1$ and $\int \xi_2(t, X(t))dt =1/2$. To avoid an identifiability problem, we may estimate $d_1$ by $\hat d_1 = n^{-1}\sum_{i=1}^n Y_i$. In the following, we will use thin-plate splines with $m=2$ to fit the data.

\subsection{Simulations}

Our first simulation study compares our estimate with other two different estimates. The first method uses the well-known functional principal component analysis (FPCA) approach. The second method uses the P-spline approach in \cite{mclean_12}, where one estimates $F$ using tensor-product B-splines with roughness penalties.
The simulation setting is the same as the setting of \cite{hall_07} and \cite{mclean_12}.
The random predictor function $X$ was generated as
$$X(t)=\zeta_1Z_1 + \sum_{k=2}^{50}\sqrt{2}~\zeta_k Z_k\cos(k\pi t),~~~ t\in [0, 1],$$ where $Z_k$ are independently sampled from the uniform distribution on $[-\sqrt{3},\sqrt{3}]$. Obviously, the $\zeta_k^2$ are eigenvalues of the covariance function of $X$. Consider two cases for the $\zeta_k$: the ''closely spaced" case and the ''well spaced" case. For the well spaced case, $\zeta_k=(-1)^{k+1}k^{-\nu/2}$ with $\nu=1.1$ and $2$. For the closely spaced case, $\zeta_1 = 1$, $\zeta_j = 0.2 (-1)^{j+1} (1 - 0.0001 j)$ for $j =2, 3, 4$, and $\zeta_{5j+k} = 0.2 (-1)^{5j+k+1}(5j)^{-{\nu}/2} -0.0001 k$ for $j\ge 1$ and $0\le k\le 4$. The true coefficient function $\beta_0$ was given by
$$\beta_0(t)=0.3 + \sum_{k=2}^{50}4\sqrt{2}(-1)^{k+1}k^{-2}\cos(k\pi t), ~~~ t\in [0, 1].$$
The simulation study was performed when the functional linear regression model is the true model. The response variable $Y$  is simulated from the model: $Y = \int_0^1 \beta_0(t)X(t)dt + \epsilon$, where the error $\epsilon\sim N(0, \sigma^2)$, where $\sigma=0.5$ and $1$. The performance of different estimators is measured by the root mean squared prediction error,
${\mathrm{RMSPE}} = \sqrt{ d^{-1} \sum_{i=1}^d\big(\widehat Y_i - Y_i\big)^2},$
where $d$ is the sample size of the test data and the $\widehat Y_i$ are predicted values.  Each training set contains $67$ curves and $33$ curves are used for the test set. For each setting, the experiment is repeated $1000$ times. The results of simulations are summarized in Table \ref{tab1}. We observe that our thin-plate spline estimator performs nearly identically to the functional PCA estimator, even though this is an ideal setting for the latter since the functional linear model holds.  Also, our estimator slightly outperforms the P-spline estimator.

\begin{table}[t]
\begin{center}
\caption{The root mean squared prediction errors (RMSPE) of three  estimators for a functional linear regression model where $Y = \int_0^1 \beta_0(t)X(t)dt + \epsilon$.
FPCA is an estimation for the functional linear model based on functional principal components analysis.  P-spline is the estimator of \cite{mclean_12}. ``ThinSpline'' is our proposed estimator using a thin-plate spline.}
\vspace{1em}
\begin{tabular}{ccc|ccc}
\hline \hline
$\xi_j$ & $\sigma$ & $\nu$  & FPCA  & P-Spline  & ThinSpline \\
\hline
\multirow{4}{*}{Well Spaced} & \multirow{2}{*}{0.5} & 1.1  & 0.61 & 0.82 &  0.68 \\
                             &                      & 2.0  & 0.52 & 0.55 &   0.56                               \\
                             &  \multirow{2}{*}{1.0}& 1.1  & 1.21 &1.65  &   1.20                                 \\
                             &                        & 2.0 & 1.04 & 1.09 &   1.08  \\ \hline
\multirow{4}{*}{Closed Spaced} & \multirow{2}{*}{0.5} & 1.1  & 0.52 & 0.53 &  0.52\\
                             &                      & 2.0  & 0.54 & 0.55 &  0.56                               \\
                             &  \multirow{2}{*}{1.0}& 1.1   &1.03 & 1.07 &   1.03                                \\
                             &                        & 2.0  & 1.06 & 1.05 &  1.04                              \\
\hline \hline
\end{tabular}
\label{tab1}
\end{center}
\end{table}

Next, we perform a simulation study to compare our estimate with the piecewise constant fit proposed in \cite{muller_12} when the additive functional model holds.
The simulation setting is the same as that in \cite{muller_12}. The predictor functions are generated according to
$$X(t) = \cos(U_1) \sin({1\over 5}\pi t) +\sin(U_1) \cos({1\over 5}\pi t)+\cos(U_2) \sin({2\over 5}\pi t)+\sin(U_2) \cos({2\over 5}\pi t)$$
for $t\in [0, 10]$ where $U_1$ and $U_2$ are iid from Uniform$[0, 2\pi]$. The sample size for the training data is $n=200$ and for the testing data is $d=1000$. The data are generated from two different nonlinear functional models: (i) $Y = \int_0^{10}\cos\{t-X(t)-5\}dt+\epsilon$; (ii) $Y=\int_0^{10} t\exp\{X(t)\}dt + \epsilon$, where $\epsilon\sim N(0, \sigma^2)$. For each setting, the experiment is repeated $50$ times. The means and the corresponding standard deviation of the root mean squared prediction error are given in Table \ref{tab2}. As expected, the functional PCA approach fails for these two examples as it has large prediction errors. In addition, our thin-plate spline estimate outperforms the piecewise constant fit (PCF) proposed in \cite{muller_12}. An additional tuning data set with sample size $200$ is used to select the needed regularization parameter in the original simulation of PCF by \cite{muller_12}. A benefit of our approach is that we do not require this tuning data set in our simulations.

\begin{table}[t]
\begin{center}
\caption{The root mean squared prediction errors (RMSPE) based on three different estimators for two nonlinear functional regression models. PCF is the piecewise constant fit of \cite{muller_12}.}
\vspace{1em}
\begin{tabular}{l|c|cccc}
\hline\hline
                                                             Model      & $\sigma $ & FPCA  &  PCF   & SSpline \\ \hline
       \multirow{3}{*}{$Y = \int_0^{10}\cos\{t-X(t)-5\}dt+\epsilon$} & 2  & 2.434 (0.018) & 2.200 (0.056) &   2.108 (0.062)     \\
                                                                      & 1  &  1.723 (0.013) & 1.156 (0.037) &  1.127 (0.035)   \\
                                                                      & 0.5 & 1.494 (0.011) & 0.680 (0.035) & 0.569 (0.026) \\
       \hline
     $Y=\int_0^{10} t\exp\{X(t)\}dt + \epsilon$                    &1 & 9.828 (0.106) & 1.119 (0.029) & 1.108 (0.031)   \\
\hline\hline
\end{tabular}
\label{tab2}
\end{center}
\end{table}

\begin{figure}[t]
\begin{center}
\resizebox{5in}{4in}{\includegraphics{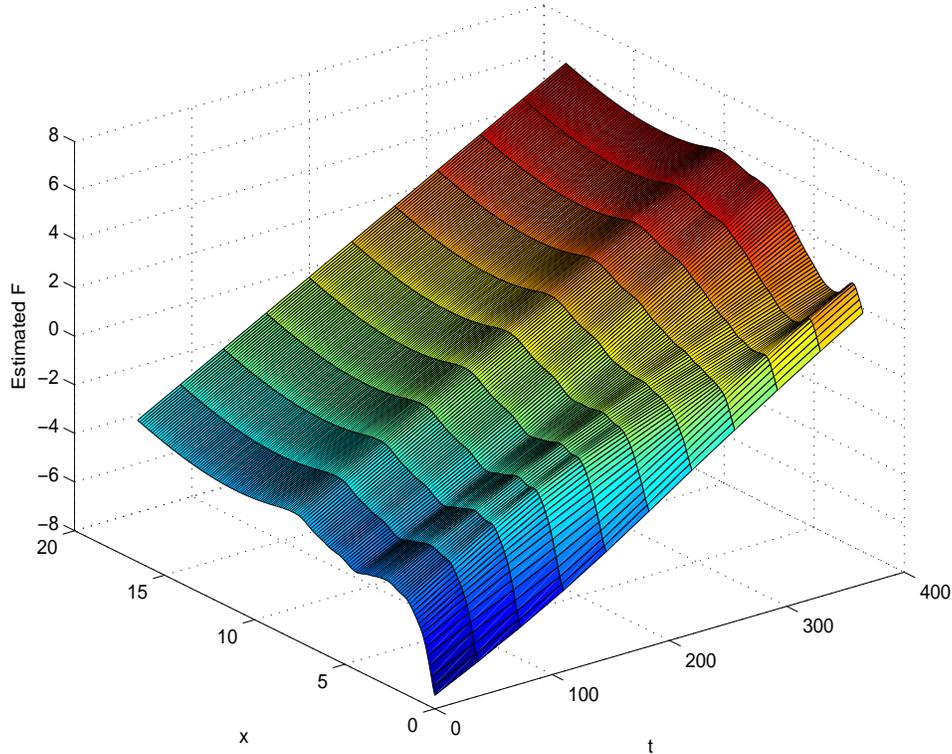}}
\end{center} \caption{Estimated surface $\widehat F_{n\lambda}(t, x)$ from the Canadian weather data.}\label{fig:weatherF}
\end{figure}

\subsection{Application: Canadian Weather Data}

The Canadian weather data example is revisited here. The dataset contains daily temperature and precipitation at 35 different locations in Canada averaged over years 1960 to 1994. Our goal is to predict the log annual precipitation based on the average daily temperature. In \cite{cai_12} it was shown that the functional PCA approach could be problematic, since the eigenfunctions corresponding to the leading eigenvalues of the covariance function seem not to represent the estimated coefficient function well. Therefore, we compare our method with the smoothing spline estimate when assuming the functional linear regression model. Under this setting, the estimate is given by
\begin{equation}\label{equ:ss}
(\hat\alpha, \hat\beta) = \arg\min \Big\{{1\over n}\sum_{i=1}^n\Big(Y_i - \alpha -\int_0^1 X_i(t)\beta(t)dt\Big)^2 + \lambda \int_0^1 (\beta''(t))^2dt\Big\}.\end{equation}

\begin{figure}[t]
\begin{center}
\resizebox{3in}{3in}{\includegraphics{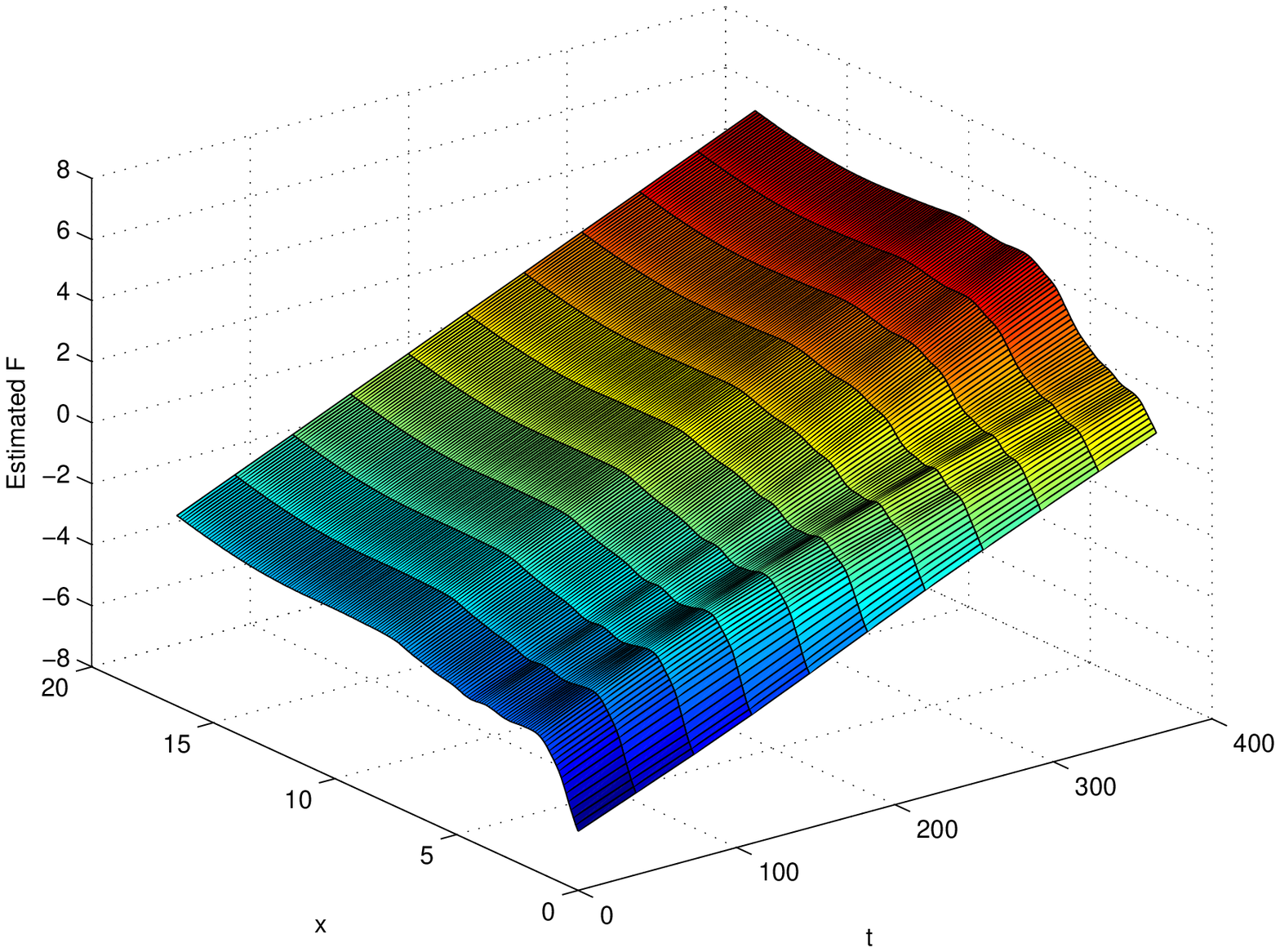}}\resizebox{3in}{3in}{\includegraphics{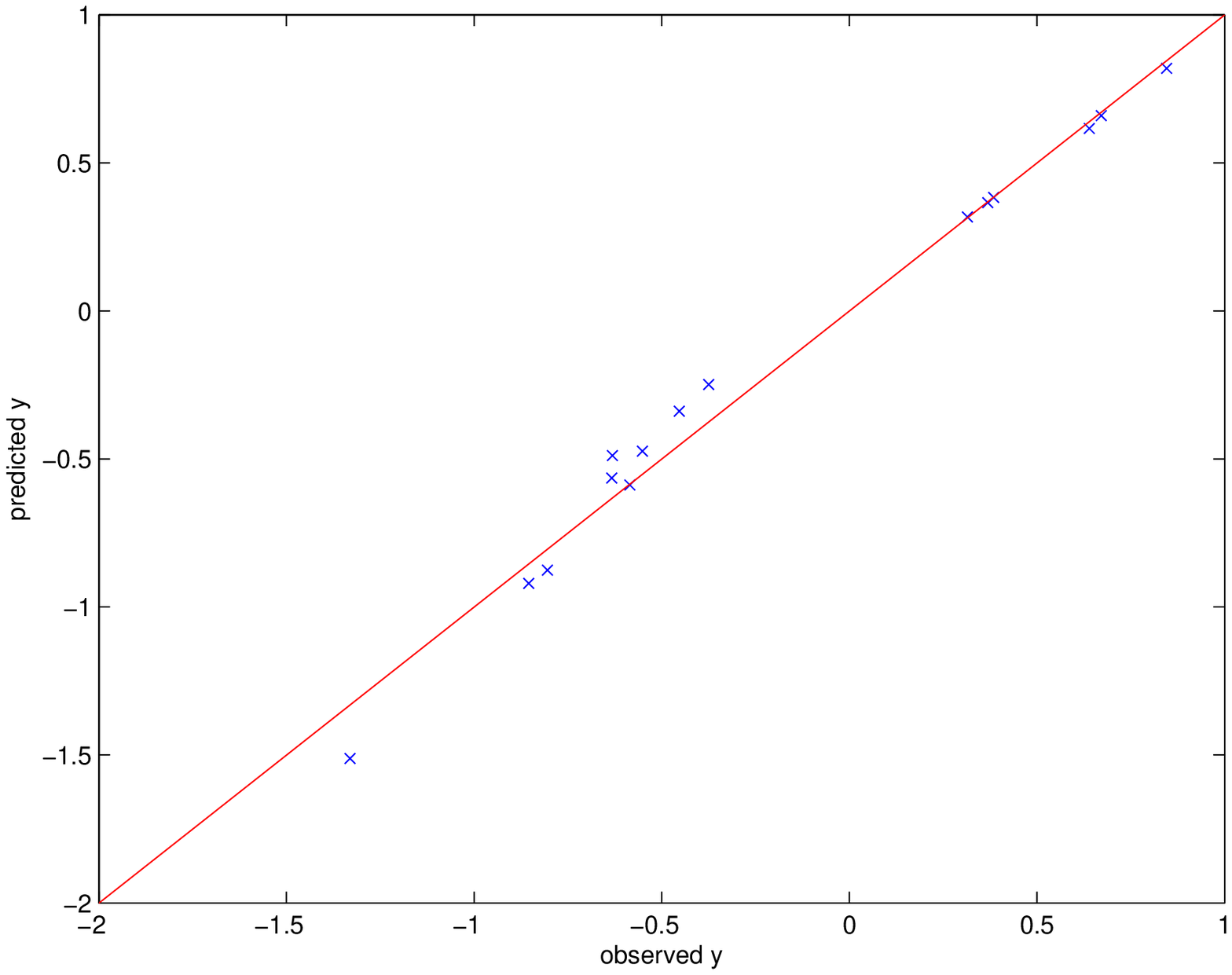}}
\end{center} \caption{Left: Estimated surface $\widehat F_\lambda(t, x)$ from the training data; Right: the predicted response versus the observed response for the testing data.}\label{fig:weatherF2}
\end{figure}

\begin{table}[t]
\begin{center}
\caption{The root mean squared prediction errors based on the estimate (\ref{equ:ss}) and the proposed predictor for Canadian weather data.}
\vspace{1em}
\begin{tabular}{c|cc}
\hline\hline
         & FLR  & ThinSpline \\ \hline
RMSPE   & 0.3014(0.1244) & 0.1110(0.0917) \\
\hline\hline
\end{tabular}
\label{tab3}
\end{center}
\end{table}

Figure \ref{fig:weatherF} shows the estimated $\widehat F_{n\lambda}$ when using the complete data.
In order to study the performance of these estimators, we randomly split the initial sample into two sub-samples: (a) A learning sample, $(X_i, Y_i)$, $i=1, \ldots, n_\ell$ with $n_\ell=20$, was used to determine the estimated coefficient function $\hat\beta_\lambda$ and the estimator $\widehat F_{n\lambda}$; (b) A test sample, $(X_i, Y_i)$, $i=n_\ell+1, \ldots, n$, with $n-n_\ell=15$ was used to evaluate the quality of the estimation. The left panel of Figure \ref{fig:weatherF2} displays the estimated $\widehat F_{n\lambda}$ from the training data set and the right panel of Figure \ref{fig:weatherF2} shows the predicted response versus the observed response for the testing data using the estimate from the training data. The points are very close to the diagonal line which indicates a good fit. We have repeated this procedure $200$ times. The mean and the corresponding standard deviations of the root mean squared prediction errors based on (\ref{equ:ss}) and our proposed predictor
are reported in Table \ref{tab3}.

It is noteworthy that the prediction error  using the continuously functional additive model is considerably less than for the functional linear regression model. The goodness-of-fit of different models is an important research topic and we will pursue this for future studies.

\subsection{Application: CA Air Quality Data}

\begin{figure}
\begin{center}
\resizebox{3in}{3in}{\includegraphics{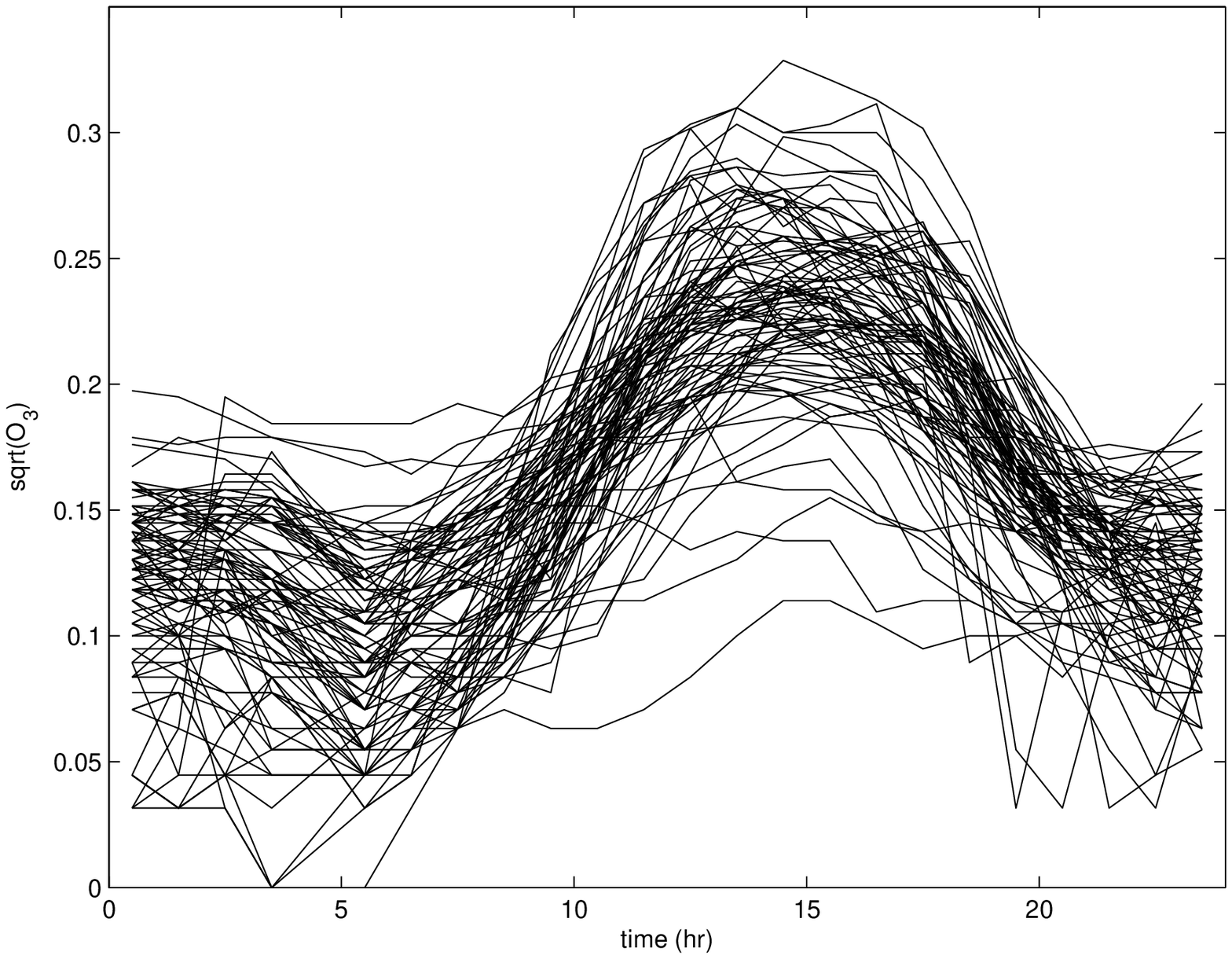}}\resizebox{3in}{3in}{\includegraphics{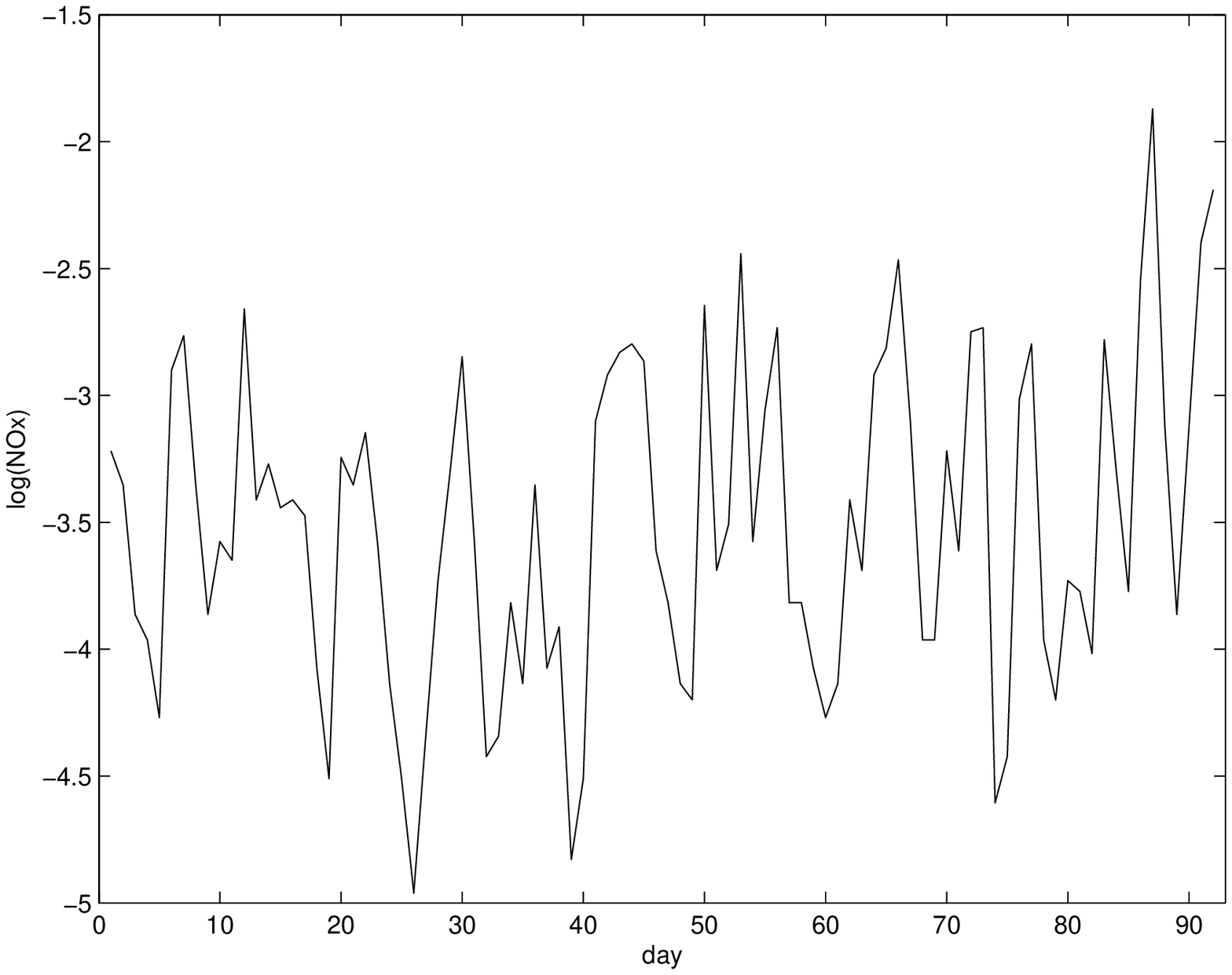}}
\end{center} \caption{Left: Daily trajectories of ground-level concentrations of ozone in the city of Sacramento in the Summer of 2005; Right: The maximum level of the ground-level concentrations of oxides of nitrogen at each day in the Summer of 2005.}\label{fig:air}
\end{figure}

Air pollutants are known to cause serious health problems. Modeling different ground level air pollutants has been an important research topics for many years. In May 2011, the California Air Resources Board has released the ``2011 Air Quality Data", which include 30 years of air quality data (1980-2009). This database, available at \url{http://www.arb.ca.gov/aqd/aqdcd/aqdcddld.htm}, contains  hourly concentrations of pollutants at different locations in California from year 1980 to year 2009. In this study, we will focus on the effect of the trajectories of ozone (O3) on the maximum level of oxides of nitrogen (NOx) in the city of Sacramento (site 3011 in the database) between June 1 and August 31 of 2005. The total sample size is $n=92$. The left panel of Figure \ref{fig:air} displays the daily trajectories of ground-level concentrations of ozone in the city of Sacramento in the Summer of 2005. For most days, we have the observations at each hour and there are a few days with some missing observations. The right panel of Figure \ref{fig:air} gives the maximum level of the ground-level concentrations of oxides of nitrogen at each day during the summer of 2005 in Sacramento.

\begin{figure}[t]
\begin{center}
\resizebox{5in}{4in}{\includegraphics{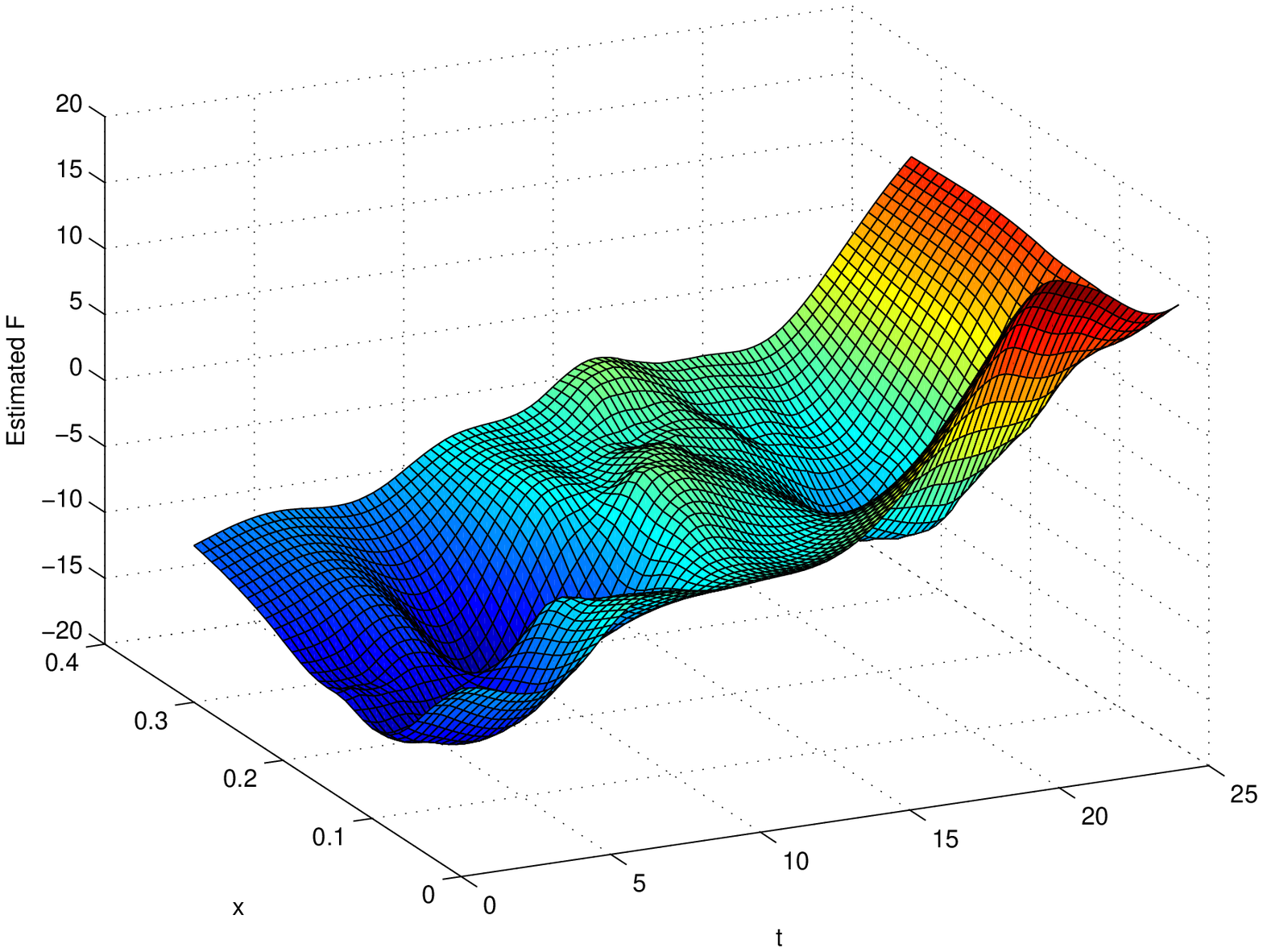}}
\end{center} \caption{Estimated surface $\widehat F_{n\lambda}(t, x)$ from the air quality data.}\label{fig:airF}
\end{figure}

\begin{figure}[t]
\begin{center}
\resizebox{3in}{3in}{\includegraphics{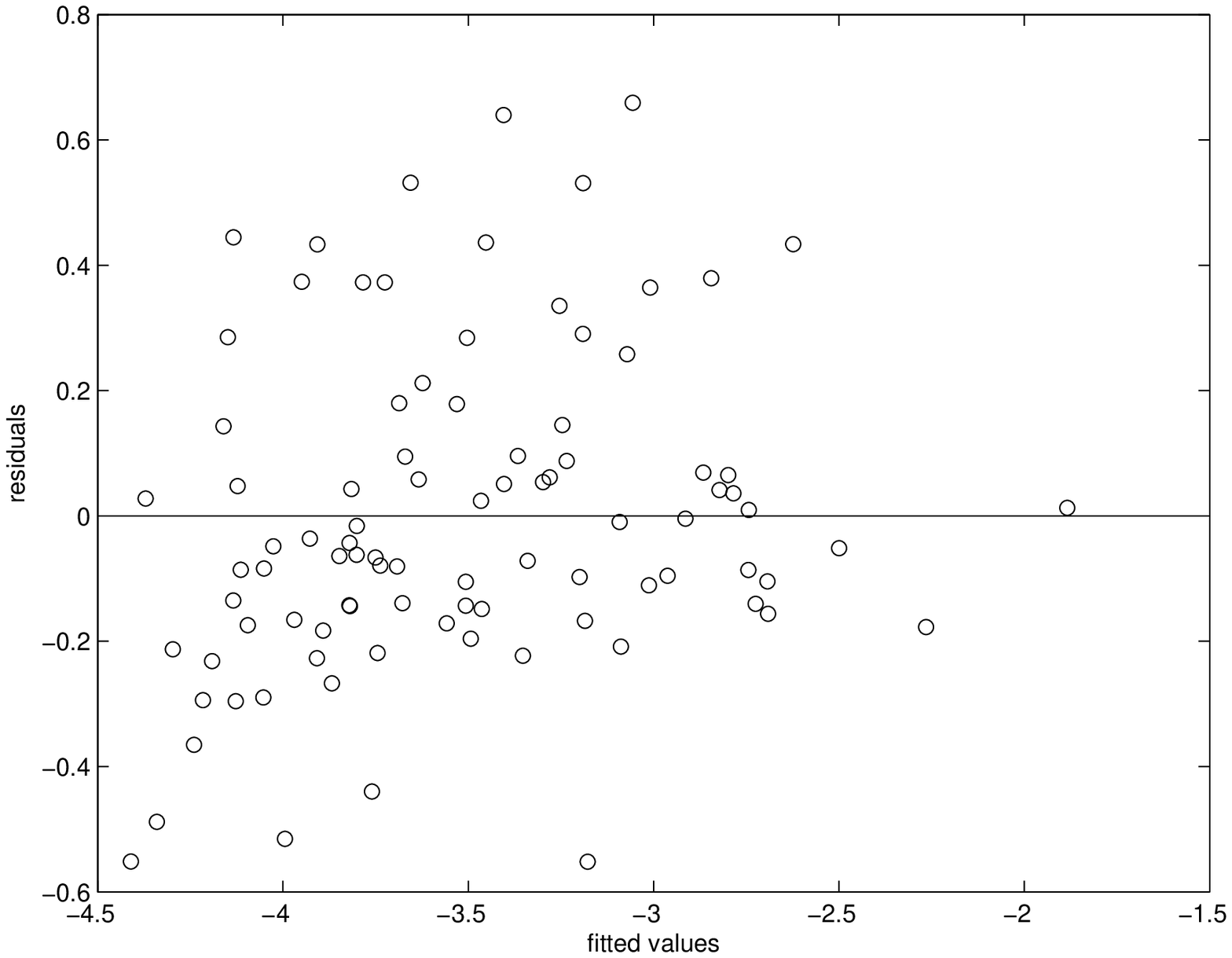}}\resizebox{3in}{3in}{\includegraphics{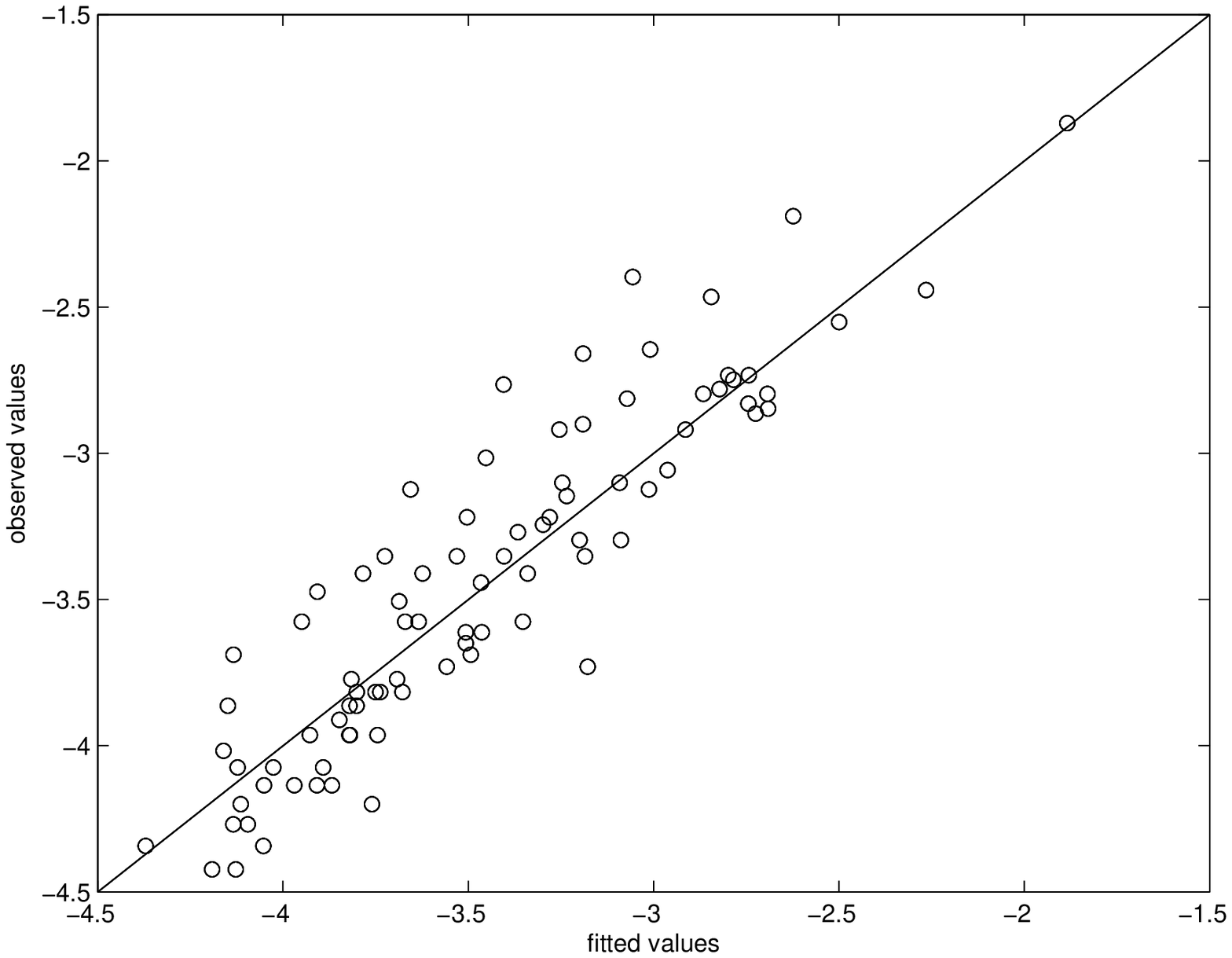}}
\end{center} \caption{Left: Residual plot; Right: Fitted values versus the observed responses.}\label{fig:airfit}
\end{figure}

Figure \ref{fig:airF} shows the estimated $\hat F_{n\lambda}$ when using the complete data. It displays a highly nonlinear pattern, which may suggest that the functional linear model may not fit the data well. To assess the goodness of fit of the additive functional model, the left panel of Figure \ref{fig:airfit} plots the residuals on the vertical axis and the fitted responses on the horizontal axis. It shows the points are randomly dispersed around the horizontal axis and did not show any typical pattern. The right panel of Figure \ref{fig:airfit} plots the fitted values versus the observed responses. The points are very closed to the diagonal line and it indicates a good fit.

We also compare the performance of the additive functional model with the functional linear regression model (\ref{equ:flr}). The $92$ observations were randomly split into training sets of size $60$ and test sets of size $32$. We repeat this procedure $1000$ times. The mean and the corresponding standard deviations of the root mean squared prediction error based on these two models are reported in Table \ref{table_air}. As expected, our additive functional linear model outperforms the functional linear model.

\begin{table}[t]
\begin{center}
\caption{The root mean squared prediction errors based on the functional linear regression (FLR) model and the additive functional model (ThinSpline) for the air quality data.}
\vspace{1em}
\begin{tabular}{c|cc}
\hline\hline
         & FLR  & ThinSpline \\ \hline
RMSPE   & 0.9450 (1.6539) & 0.6148(0.0985) \\
\hline\hline
\end{tabular}
\label{table_air}
\end{center}
\end{table}

\section{Discussion}

We have established the minimax rate of convergence for prediction for the continuous functional additive model. It is shown that the optimal rate depends on the decay rate of the eigenvalues of the kernel $C$, which depends on the reproducing kernel and the joint distribution of the random predictor function  at any two points. The minimax theory in the existing literature on the functional linear regression model is a special setting of current study.

We have focused on the additive functional  model with the squared error loss in this paper. It should be noted that the method of regularization can be easily extended to handle other models such as the generalized regression model \cite{cardot_05, muller_05,mclean_12, du_12}. We shall leave these extensions for future papers.

The simulation  in this paper study only the estimator using thin-plate splines. For the case of univariate regression, \cite{wang_11} has showed that a smoothing spline and a P-spline are asymptotically equivalent. Similar asymptotic equivalent result is expected to hold for the bivariate regression too. So, it is expected that our simulation performance is similar to that of \cite{mclean_12}, who used the bivariate P-splines to fit the data. However, it should be pointed out that our results can be applied to the more general reproducing kernel Hilbert spaces.

It is worth noting that estimating $F_0$ itself is totally different problem with the prediction discussed in the current paper. For example, for the functional linear regression model, we may not estimate the coefficient function $\beta_0$ consistently without additional conditions linking the smoothness of $\beta_0$ and the curves $X_i$ \cite{crambes_09}. As an example of additional assumptions, one might assume the reproducing kernel $K$ and the covariance kernel $G$ are perfectly aligned, i.e., they share the same set of eigenfunctions. Under this circumstance, we may estimate $\beta_0$ consistently \cite{yuan_10}. It deserves further study when we can estimate $F_0$ consistently under the additive functional  model.
This issue is important and we could use this to test for linearity of $F_0$.

\section{Proofs}

\subsection{Proof of Theorem \ref{thm:low}}

In the following proofs, let $c_i$, $i=1,2,\ldots$ be generic constants which change from line to line.

Since any lower bound for a specific case yields immediately a lower bound for the general case, to establish lower bounds, we only study the case when the $\epsilon_i$ are i.i.d.\ $N(0,\sigma^2)$. Fix $\alpha\in (0, 1/8)$. It follows from Theorem 2.5 in \cite{tsybakov_09} that in order to establish the minimax lower bound for ${\mathfrak R}_n$, for each $n$ we need to find
functions $\{F_{jn}$, $j=0,\ldots,M\}$, satisfying   the following three conditions:
\begin{enumerate}
\item[(a).] $F_{jn} \in {\cal H}(K)$, $j=0, \ldots, M$,
\item[(b).] $\mathbb E^*\Big\{ \int_0^1 \Big[  F_{jn}(t, X_{n+1}(t)) - F_{kn}(t, X_{n+1}(t))   \Big]dt \Big\}^2\ge 2s$, \newline for $0\le j<k\le M$,
    \item[(c).] $M^{-1}\sum_{j=1}^M{\cal K}(P_j, P_0)\le \alpha \log M$, where $P_j$ denotes the joint distribution of $\{(Y_i, X_i):i=1,\ldots n\}$ when $F_0 = F_{jn}$ and  ${\cal K}(\cdot, \cdot)$ is the Kullback-Leibler distance between two probability measures.
\end{enumerate}
We will specify $M \to \infty$ and $s \to 0$ later. If (a), (b), and (c) are satisfied, then the minimax lower bound for the rate of convergence of  ${\mathfrak R}_n$ has the same order as $s$.

First we verify part (a). Let $m$ be the smallest integer greater than $c_0n^{1/(2r+1)}$ for some positive constant $c_0$ to be specific later.
For a $\omega = (\omega_{m+1}, \ldots, \omega_{2m}) \in \{0, 1\}^m$, let
$$F_\omega = \sum_{j=m+1}^{2m} \omega_j m^{-1/2} K^{1/2}(\phi_j).$$


$F_\omega \in {\cal H}(K)$ for all $\omega$ if $K^{1/2}(\phi_j) \in {\cal H}(K)$ for all $j$. Thus, we need to show that $\big\langle K^{1/2}(\phi_j),K\big(\cdot,(t,x)\big)\big\rangle = K^{1/2}(\phi_j)(t,x)$. This result holds since
\[
\big\langle K^{1/2}(\phi_j),K\big(\cdot,(t,x)\big)\big\rangle
= \big\langle K(\phi_j),K^{1/2}\big(\cdot,(t,x)\big)\big\rangle
= \big\langle \phi_j,K^{1/2}\big(\cdot,(t,x)\big)\big\rangle
= K^{1/2}(\phi_j)(t,x).
\]
We also have
\begin{align*}
\langle{K^{1/2}}(\phi_j), {K^{1/2}}(\phi_k)\rangle_{{\cal H}(K)} = \langle \phi_j, K(\phi_k)\rangle_{{\cal H}(K)}
=\langle\phi_j, \phi_k\rangle_{L_2} = \delta_{jk},
\end{align*}
where $\delta_{jk}=1$ for $j=k$, and $0$ for $j\neq k$.

Further, the Varshamov-Gilbert bound
(see \cite{tsybakov_09}, p.\ 104)
shows that, for $m\ge 8$, there exists a subset $\Omega =\{\omega^{0}, \omega^{1}, \ldots, \omega^{M}\}\subseteq \{0, 1\}^m$ such that $\omega^{0}=\{0,  \ldots, 0\}$,
\begin{equation}
d(\omega^{j}, \omega^{k})\ge {m\over 8}, ~~~\forall~ 0\le j<k\le M,\label{dr1}
\end{equation}
where  $d(\omega^j, \omega^k)=\sum_{i=m+1}^{2m} I(\omega^j_i\ne \omega^k_i)$ is the Hamming distance between $\omega^{j}$ and $\omega^{k}$, and
$$M \ge 2^{m/8}.$$

To verify part (b), for $\omega, \omega'\in \Omega$, direct calculation yields that
\begin{align*}
&~\mathbb E^*\Big\{ \int_0^1 \Big[  F_{\omega}(t, X_{n+1}(t)) - F_{\omega'}(t, X_{n+1}(t))   \Big]dt \Big\}^2\\
= & \sum_{j=m+1}^{2m}\sum_{k=m+1}^{2m}  m^{-1} (\omega_j - \omega_j')(\omega_k - \omega_k')\int\int \mathbb E^* \Big[ {K^{1/2}}(\phi_j)(t, X(t)) {K^{1/2}}(\phi_k)(s, X(s))\Big]dtds\\
=& \sum_{k=m+1}^{2m}  m^{-1} (\omega_k - \omega_k')^2 \rho_k
\ge ~  m^{-1}\rho_{2m} d(\omega, \omega') \ge ~ c_1 m^{-1} (2m)^{-2r} m/8 \ge c_2 n^{-2r/(2r+1)}
\end{align*}
by \eqref{dr1}, $\rho_k \asymp k^{-2r}$, and the definition of $m$.
Hence, $s$ in part (b) is of order $n^{-2r/(2r+1)}$.

Next, observe that for any $\omega, \omega' \in \Omega$,
\begin{align*}
\log(P_{F_{\omega'}}/P_{F_\omega}) = {1\over \sigma^2} \sum_{i=1}^n &\Big( Y_i - \int F_\omega(t, X(t))dt\Big)\int \Big\{F_\omega(t, X(t)) - F_{\omega'}(t, X(t))\Big\}dt -\\
&~~~~ {1\over 2\sigma^2}
\sum_{i=1}^n\Big[ \int\Big\{F_\omega(t, X(t)) - F_{\omega'}(t, X(t))\Big\}dt  \Big]^2.
\end{align*}
Therefore,
\begin{align*}
& \qquad  {\cal K}(P_{F_{\omega'}}, P_{F_{\omega}}) = {n\over 2\sigma^2} \mathbb E^* \Big[ \int\Big\{F_\omega(t, X(t)) - F_{\omega'}(t, X(t))\Big\}dt  \Big]^2\\
& = {n\over 2\sigma^2} \sum_{k=m+1}^{2m}m^{-1} (\omega_k - \omega_k')^2 \rho_k \le {n\over 2\sigma^2} \rho_m \sum_{k=m+1}^{2m}m^{-1} (\omega_k - \omega_k')^2 \le {n\over 2\sigma^2}m^{-2r} \le c_3 n^{1/(2r+1)}.
\end{align*}
Since $m$ is the smallest integer greater than $c_0n^{1/(2r+1)}$, this implies that
$${1\over M}\sum_{j=1}^M{\cal K}(P_j, P_0) \le c_3 n^{1/(2r+1)} \le \alpha \log M,$$
if we choose $c_0\ge 8c_3/(\alpha\log 2)$ and $M = 2^{m/8}$.
This completes the proof of Theorem \ref{thm:low}.

\subsection{Proofs of Theorem \ref{thm:sol} and Theorem \ref{thm:upper}}

\noindent {\it Proof of Theorem \ref{thm:sol}}.
Define the subspace of ${\cal H}$,
$$\overline{\cal H}_1 = \mbox{span}\Big\{ \int K\Big((t, x); (s, X_i(s))\Big)ds: i=1,\ldots, n\Big\}.$$
Note that $\overline{\cal H}_1$ is a closed linear subspace of ${\cal H}_1$. For any $F\in {\cal H}$, one may write
$$F = F_0+F_1+\delta,$$
where $F_0\in {\cal H}_0$, $F_1\in \overline{\cal H}_1$ and  $\delta\in {\cal H}_1 \ominus \overline{\cal H}_1$.
%
Observe that
$$\eta_F(X_i) = \int F(t, X_i(t))dt = \eta_{F_0+F_1}(X_i),$$
because
$$\eta_\delta(X_i) = \Big\langle \int K\Big((\cdot; (s, X_i(s))\Big)ds, \delta\Big\rangle_{\cal H} = 0.$$
Further, due to orthogonality, $\|F\|_{\cal H}^2 = \|F_0+F_1\|_{\cal H}^2+\|\delta\|_{\cal H}^2$ and
$\|F_0+F_1\|_{\cal H}^2 \le \|F\|_{\cal H}^2$. Therefore, the minimum of (\ref{equ:obj}) must belong to the linear space ${\cal H}_0\oplus \overline{\cal H}_1$. \qed \\

\noindent  {\it Proof of Theorem \ref{thm:upper}}. Note that $L_2(K^{1/2}) = {\cal H}(K)$. So there exist $G_0$ and $\hat G_\lambda$ such that $F_0 = K^{1/2}(G_0)$ and $\widehat F_{n\lambda} = K^{1/2}(\hat G_\lambda)$. Therefore,
\begin{align*}
\eta_{F_0}(X) &= \int F_0(t, X(t))dt = \int \Big\langle K\Big(\cdot; (s, X(s)\Big), F_0\Big\rangle_{\cal H(K)}ds\\
&= \int \Big\langle K^{1/2}\Big(\cdot; (s, X(s))\Big), G_0 \Big\rangle_{L_2}ds,
\end{align*}
and
$$\mathfrak{R}_n = \mathbb E^* \Big| \int \Big\langle K^{1/2}\Big(\cdot; (s, X(s))\Big), \hat G_\lambda -G_0 \Big\rangle_{L_2}ds  \Big|^2= \Big\|\hat G_\lambda - G_0\Big\|_C^2,$$
where
$$
\Big\|G\Big\|_C^2 = \int\cdots\int G\Big((t, x);(u_1, z_1)\Big)  C\Big((u_1, z_1); (u_2, z_2)\Big)  G\Big((u_2, z_2);(s, y)\Big).
$$
Write
$$C_n\Big((t, x); (s, y)\Big) = {1\over n}\sum_{i=1}^n \int\int K^{1/2}\Big( (t, x);(u, X_i(u))  \Big)K^{1/2}\Big( (s, y);(v, X_i(v))  \Big)dudv.$$
Recall that $Y_i = \int \Big\langle K^{1/2}\Big(\cdot; (s, X_i(s))\Big), G_0\Big\rangle ds + \epsilon_i$. Denote
$g_n = {1\over n}\sum_{i=1}^n \epsilon_i \int K^{1/2}\Big(\cdot; $ $(s, X(s))\Big)ds.$
Then,
$\hat G_\lambda = \Big( C_n + \lambda I \Big)^{-1} \Big( C_n(G_0) + g_n  \Big). $
Define
$G_\lambda = \Big( C + \lambda I \Big)^{-1} C(G_0).$
It follows from triangle inequality that
\begin{equation}\label{equ:G}
\Big\|\hat G_\lambda - G_0\Big\|_C \le \Big\|G_\lambda - G_0\Big\|_C + \Big\|\hat G_\lambda - G_\lambda\Big\|_C.
\end{equation}
Let us first bound the first term in the right hand side of (\ref{equ:G}). Recall that the $\phi_k$ are the eigenfunctions of $C$. Write $G_0 = \sum_{k=1}^\infty a_k\phi_k$. Then,
$$G_\lambda = \sum_{k=1}^\infty {a_k\rho_k\over \lambda + \rho_k}\phi_k,$$
and
$$\Big\|G_\lambda - G_0\Big\|_C^2 =  \sum_{k=1}^\infty {\lambda^2a_k^2\rho_k\over (\lambda + \rho_k)^2}\le \lambda^2 \max_{k\ge 1}{\rho_k\over (\lambda + \rho_k)^2}\sum_{k=1}^\infty a_k^2 = O(\lambda) \Big\|G_0\Big\|_{L_2}^2.$$

Next, let us bound the second term in the right hand side of (\ref{equ:G}). Recall that
$\Big( C_n + \lambda I \Big) \hat G_\lambda = C_n(G_0) + g_n.$
We observe that
\begin{align*}
G_\lambda - \hat G_\lambda &= (C+ \lambda I)^{-1} (C_n + \lambda I)(G_\lambda - \hat G_\lambda) + (C+ \lambda I)^{-1} (C - C_n )(G_\lambda - \hat G_\lambda)\\
&= (C + \lambda I)^{-1} (C_n+\lambda I) G_\lambda -  (C + \lambda I)^{-1} C_n G_0 - (C + \lambda I)^{-1}  g_n \\
&~~~~~~~~~~~~~~~~~~+(C+ \lambda I)^{-1} (C - C_n )(G_\lambda - \hat G_\lambda)\\
&=(C + \lambda I)^{-1} C_n ( G_\lambda -G_0) + \lambda (C+ \lambda I)^{-2} CG_0 - (C + \lambda I)^{-1} g_n \\
&~~~~~~~~~~~~~~~~~~+(C+ \lambda I)^{-1} (C - C_n )(G_\lambda - \hat G_\lambda)\\
& = (C + \lambda I)^{-1} C ( G_\lambda -G_0) + \lambda (C+ \lambda I)^{-2} CG_0 - (C + \lambda I)^{-1} g_n \\
&~~~~~~~~~~~~~~~~~~+(C + \lambda I)^{-1} (C_n-C) ( G_\lambda -G_0)\\
&~~~~~~~~~~~~~~~~~~+(C+ \lambda I)^{-1} (C - C_n )(G_\lambda - \hat G_\lambda)
= \text{I + II + III + IV + V}.\\
\end{align*}
We now bound five terms on the right hand side separately. Direct calculation yields that
\begin{align*} &
\big\|\text{I}\big\|^2_C=
\Big\|(C + \lambda I)^{-1} C ( G_\lambda -G_0)\Big\|_C^2  = \lambda^2 \sum_{k=1}^\infty {a_k^2\rho_k^3\over (\lambda+ \rho_k)^4}\\
& \le \lambda^2 \max_{k\ge 1}{\rho_k^3\over (\lambda+ \rho_k)^4} \sum_{k=1}^\infty a_k^2 = O(\lambda) \Big\|G_0\Big\|_{L_2}^2.
\end{align*}
Similarly,
$$\big\|\text{II}\big\|^2_C=\Big\|\lambda (C+ \lambda I)^{-2} CG_0\Big\|_C^2 = \lambda^2 \sum_{k=1}^\infty {a_k^2\rho_k^3\over (\lambda+ \rho_k)^4} \le O(\lambda) \Big\|G_0\Big\|_{L_2}^2.$$

Next, we make use three auxiliary results whose proofs are similar to ones in Cai and Yuan (2012) so we omit the details. If there exists a constant $c>0$ such that
$$\mathbb E\Big( \int F(t, X(t))dt \Big)^4 \le c\Big( \mathbb E\Big( \int F(t, X(t))dt \Big)^2  \Big)^2,$$
for any $\nu>0$ such that $2r (1-2\nu)>1$, then
\begin{equation}\label{equ:use1}\Big\|C^{\nu}(C+ \lambda I)^{-1} (C - C_n ) C^{-\nu}\Big\|_{\mathrm op} = O_p\Big( \Big(n\lambda^{1-2\nu+1/(2r)}\Big)^{-1/2}  \Big),\end{equation}
and
\begin{equation}\label{equ:use2}\Big\|C^{1/2}(C+ \lambda I)^{-1} (C - C_n ) C^{-\nu}\Big\|_{\mathrm op} = O_p\Big( \Big(n\lambda^{1/(2r)}\Big)^{-1/2}  \Big),\end{equation}
where $\|\cdot\|_{\mathrm op}$ stands for the usual operator norm. Further, for any $0\le \nu\le 1/2$
\begin{equation}\label{equ:use3}\Big\| C^{\nu}(C+\lambda I)^{-1} g_n\Big\|_{L_2} = O_p\Big(\Big( n\lambda^{1-2\nu+1/(2r)}  \Big)^{-1/2}\Big).
\end{equation}
Using (\ref{equ:use1}) we have
\begin{align*}
\Big\|C^{\nu}(C+ \lambda I)^{-1} (C - C_n )(G_\lambda - \hat G_\lambda)\Big\|_{L_2}^2 & \le \Big\|C^{\nu}(C+ \lambda I)^{-1} (C - C_n )C^{-\nu} \Big\|_{\mathrm op} ~\Big\|C^{\nu}(G_\lambda - \hat G_\lambda)\Big\|_{L_2}^2\\
& \le o_p(1) \Big\|C^{\nu}(G_\lambda - \hat G_\lambda)\Big\|_{L_2}^2,
\end{align*}
whenever $\lambda \ge cn^{-2r/(2r+1)}$ for some constant $c>0$. Similarly,
\begin{align*}
\Big\|C^{\nu}(C+ \lambda I)^{-1} (C - C_n )(G_\lambda - G_0)\Big\|_{L_2}^2 & \le \Big\|C^{\nu}(C+ \lambda I)^{-1} (C - C_n )C^{-\nu} \Big\|_{\mathrm op} ~\Big\|C^{\nu}(G_\lambda - G_0)\Big\|_{L_2}^2\\
& \le o_p(1) \Big\|C^{\nu}(G_\lambda - G_0)\Big\|_{L_2}^2.
\end{align*}
So, for $0<\nu<1/2-1/(4r)$,
\begin{align*}
\Big\|& C^{\nu} (G_\lambda - \hat G_\lambda)  \Big\|_{L_2} \le \Big\|C^{\nu}(C+ \lambda I)^{-1} C(G_\lambda - G_0)\Big\|_{L_2}+ \Big\|C^{\nu}(C+ \lambda I)^{-1} (C - C_n )(G_\lambda - G_0)\Big\|_{L_2}\\
&+ \lambda \|C^{1+\nu}G_0\|_{L_2} + \|C^{\nu}(C+ \lambda I)^{-1}g_n\|_{L_2}+\Big\|C^{\nu}(C+ \lambda I)^{-1} (C - C_n )(G_\lambda - \hat G_\lambda)\Big\|_{L_2}\\
& = O_p\Big( \lambda^\nu + \Big( n\lambda^{1-2\nu+1/(2r)}  \Big)^{-1/2} \Big) = O_p(\lambda^\nu),
\end{align*}
when $c_1n^{-2r/(1+2r)} \le \lambda \le c_2n^{-2r/(1+2r)}$ for $0<c_1<c_2<\infty$.
Next,
\begin{align*}
& \big\|\text{IV}\big\|_C=
\Big\|(C + \lambda I)^{-1} (C_n-C) ( G_\lambda -G_0)\Big\|_C  =  \Big\|C^{1/2}(C + \lambda I)^{-1} (C_n-C) ( G_\lambda -G_0)\Big\|_{L_2}\\
&  \qquad\le  \Big\|C^{1/2}(C + \lambda I)^{-1} (C_n-C) C^{-\nu}\Big\| \|T^\nu (G_\lambda -G_0)\|_{L_2} \\ &\qquad \le O_p\Big((n\lambda^{1/(2r)})^{-1/2}\lambda^\nu \Big) = o_p\Big( (n\lambda^{1/(2r)})^{-1/2} \Big).
\end{align*}
Similarly,
\begin{align*}
&
\big\|\text{V}\big\|_C=
\Big\|(C + \lambda I)^{-1} (C_n-C) (G_\lambda - \hat G_\lambda)\Big\|_C  =  \Big\|C^{1/2}(C + \lambda I)^{-1} (C_n-C) ( G_\lambda - \hat G_\lambda)\Big\|_{L_2}\\
&\le  \Big\|C^{1/2}(C + \lambda I)^{-1} (C_n-C) C^{-\nu}\Big\| \|T^\nu (G_\lambda - \hat G_\lambda)\|_{L_2}  \le O_p\Big((n\lambda^{1/(2r)})^{-1/2}\lambda^\nu \Big)\\ & \qquad = o_p\Big( (n\lambda^{1/(2r)})^{-1/2} \Big).
\end{align*}
It follows from (\ref{equ:use3}),
$$\big\|\text{III}\big\|_C=
\Big\|(C + \lambda I)^{-1} g_n\Big\|_C  = \Big\|C^{1/2}(C + \lambda I)^{-1} g_n \Big\|_{L_2} = O_p\Big((n\lambda^{1/(2r)})^{-1/2} \Big).$$
Combining the facts above, we conclude that, if $\lambda$ is of order  $n^{-{2r\over 2r+1}}$, then
$\| G_\lambda - \hat  G_\lambda \|_C = O_P(n^{-{2r\over 2r+1}}).$
\qed

\end{document}